\documentclass[conference, letterpaper]{IEEEtran}
\IEEEoverridecommandlockouts
\usepackage{cite}
\usepackage{amsmath,amssymb,amsfonts}
\usepackage{algorithmic}
\usepackage{graphicx}
\usepackage{textcomp}
\usepackage{xcolor}
\def\BibTeX{{\rm B\kern-.05em{\sc i\kern-.025em b}\kern-.08em
    T\kern-.1667em\lower.7ex\hbox{E}\kern-.125emX}}
    
\usepackage{tikz}
\usepackage{pgfplots} 
\usetikzlibrary{external}
\DeclareMathOperator*{\argmin}{argmin}
\def\tr{\mathrm{tr}}

\newtheorem{assumption}{Assumption}
\newtheorem{lemma}{Lemma}

\begin{document}

\title{\LARGE \bf
Imitation Learning of Stabilizing Policies for Nonlinear Systems
}


\author{Sebastian East
\thanks{S.~East is with the Department of Aerospace Engineering, University of Bristol, BS8 1TR, {\tt\footnotesize sebastian.east@bristol.ac.uk}.}%
}

\maketitle

\begin{abstract}
There has been a recent interest in imitation learning methods that are guaranteed to produce a stabilizing control law with respect to a known system. Work in this area has generally considered linear systems and controllers, for which stabilizing imitation learning takes the form of a biconvex optimization problem. In this paper it is demonstrated that the same methods developed for linear systems and controllers can be readily extended to polynomial systems and controllers using sum of squares techniques. A projected gradient descent algorithm and an alternating direction method of multipliers algorithm are proposed as heuristics for solving the stabilizing imitation learning problem, and their performance is illustrated through numerical experiments.
\end{abstract}

\section{Introduction}

Imitation learning is the process of synthesising a controller to approximate the behaviour of an `expert' that can successfully complete a task. This approach is often taken for control problems that can be readily solved by a human, but pose significant challenges when using traditional controller design techniques (e.g. operating a vehicle). In principle, imitation learning is a supervised learning problem, where the objective is to determine a function that maps a set of known states to a corresponding set of known control inputs, so powerful function approximation techniques can readily be applied. In practice, however, ensuring that a system behaves satisfactorily when controlled by a learned policy is not trivial, and is an open area of research \cite{Hussein2017}.

Recently, there has been interest in developing imitation learning methods that can guarantee that the learned policy is \emph{stabilizing} when used to control a system with known dynamics; a problem that has obvious value when considering controller design for safety-critical applications. Palan et. al. \cite{Palan2020} posed the imitation learning problem with a `Kalman constraint' that ensures the learned controller is optimal for a known linear time-invariant (LTI) system and unknown quadratic cost function, and is therefore stabilizing (this approach has strong parallels with \emph{inverse optimal control}: the problem of identifying the cost function used to generate a control law for a given system, e.g. \cite{Menner2018}). The imitation learning problem in \cite{Palan2020} was biconvex, and the alternating direction method of multipliers (ADMM) \cite{ADMM} was proposed as a heuristic solution algorithm. Havens and Hu \cite{Havens2021} relaxed the (potentially conservative) Kalman constraint and instead used linear matrix inequalities to ensure stability and robustness for a known LTI system. The resulting imitation learning problem was also biconvex, and a projected gradient descent algorithm \cite[\S 3.1]{Bubeck2015} was investigated as an alternative to ADMM. Yin et. al. \cite{Yin2022} used a similar approach to ensure that a neural network controller trained with imitation learning was stabilizing w.r.t. a LTI system, but the ADMM algorithm proposed for that approach required the repeated solution of a nonconvex optimization problem, which can itself only be solved locally. Imitation learning can also be approached using \emph{differentiable control} (where an optimization-based control policy is differentiated so that its parameters can be learned \cite{Amos2018DifferentiableMF}), and in \cite{East2020} it was demonstrated that a constrained linear quadratic control policy can be analytically differentiated, and therefore trained for imitation learning using gradient-based optimization.

A major limitation of the previous work in this area is that stability has only been considered with respect to LTI systems, and often for linear state-feedback controllers. Systems typically only behave linearly within a neighbourhood of a chosen set-point, so the stability certificates obtained using the aforementioned approaches may only hold in a small region of state space. Furthermore, the restriction to linear control laws is conservative, as a nonlinear controller may be required to stabilize a nonlinear system, and many well known synthetic control methodologies are nonlinear (e.g. the constrained linear quadratic regulator is piecewise affine \cite{BEMPORAD20023}). In this paper this limitation is therefore addressed with an imitation learning method for nonlinear controllers that are guaranteed to stabilize a known nonlinear system. In particular, it is demonstrated that the approach considered in \cite{Havens2021} can be readily extended to \emph{polynomial} systems and controllers using sum of squares techniques developed for optimal control synthesis \cite{Prajna2004}. The ADMM and projected gradient algorithms proposed in \cite{Havens2021} are reevaluated here in numerical experiments, where it is found that ADMM generally has superior performance and that, in the case of a linear controller, ADMM can obtain a highly accurate solution in a handful of iterations. The approach here is similar to that recently presented in \cite{Guo2020}, but that paper demonstrated a method for synthesising a controller from a finite number of state measurements of a polynomial system, and did not consider imitation learning.

\subsection{Preliminaries}
The natural numbers, integers, and real numbers are given by $\mathbb{N}$, $\mathbb{Z}$, and $\mathbb{R}$. $\{ a, \dots, b \}$ denotes the set of integers between $a \in \mathbb{Z}$ and $b \in \mathbb{Z}$ (inclusive), and $\mathbb{S}^{n}$ denotes the real $n \times n$ symmetric matrices. The notation $\otimes$ represents the Kronecker product. For $x \in \mathbb{R}^n$ and $\alpha \in \mathbb{N}^n$, the monomial $x^{\alpha} := x_1^{\alpha_1} x_2^{\alpha_2} \dots x_n^{\alpha_n}$ has degree $|\alpha |:= \sum_{i=1}^n \alpha_i$. The set of all monomials of $x \in \mathbb{R}^n$ with degree $\leq d$ is denoted $\mathcal{M}_d[x]$ and has cardinality $|\mathcal{M}_d[x]| = \sum_{i=0}^d \begin{pmatrix} i + n - 1 \\ n - 1 \end{pmatrix}$. The notation $\mathcal{M}_d[x]_i$ represents the $i$-th element of $\mathcal{M}_d[x]$ determined using an arbitrary (but fixed) ordering. A function $f: \mathbb{R}^n \mapsto \mathbb{R}$ is a polynomial iff $f(x) = \sum_{i=1}^N c_i x^{\alpha(i)}$ where $c_i \in \mathbb{R}$ and $\alpha(i) \in \mathbb{N}^n$ for all $i \in \{1, \dots, N \}$, and the degree of $f$ is $\max_{i \in \{1, \dots, N \}} |\alpha(i) | $. The set of all polynomials is denoted $\mathcal{P}$, and the set of polynomials with degree less than or equal to $d$ is denoted $\mathcal{P}_{d}$. A vector-valued function $v: \mathbb{R}^n \mapsto \mathbb{R}^m$ is polynomial iff $v_i \in \mathcal{P}$ for all $i \in \{1, \dots, m\}$, and the degree of $v$ is the maximum degree of $v_i$. The set of all vector-valued polynomials with range $\mathbb{R}^m$ is denoted $\mathcal{P}^m$, and the elements of $\mathcal{P}^m$ with degree less than or equal to $d$ is denoted $\mathcal{P}^m_d$. A matrix function $M : \mathbb{R}^n \mapsto \mathbb{R}^{m \times l}$ is polynomial if $M_{i,j} \in \mathcal{P}$ for all $i \in \{1, \dots, m \}$ and $j \in \{1, \dots, l \}$, and the degree of $M$ is the maximum degree of $M_{i,j}$. The set of all matrix polynomials with range $\mathbb{R}^{m \times n}$ is denoted $\mathcal{P}^{m \times n}$, and the elements of $\mathcal{P}^{m \times n}$ with degree less than or equal to $d$ is denoted $\mathcal{P}^{m \times n}_d$. Any matrix polynomial $M \in \mathcal{P}^{m \times n}_d$ has the equivalent representation $ M(x) = \sum_{i=1}^{|\mathcal{M}_d[x]|} M_i \mathcal{M}_d[x]_i$, where $M_i \in \mathbb{R}^{m \times n}$ for all $i$.

\section{Problem Formulation}

\subsection{Imitation Learning}

In this paper, imitation learning is framed as a supervised learning problem. Consider a system described by a state $x \in \mathbb{R}^n$ and control input $u \in \mathbb{R}^m$, and assume that there exists $N$ sampled pairs of input-output data, $(\hat{x}, \hat{u})$, generated by an expert controlling the system:
\begin{equation}\label{eqn::data}
\{(\hat{x}_1, \hat{u}_1 ), \dots, (\hat{x}_N, \hat{u}_N)\}.
\end{equation}
This data need not be ordered, and could be generated from individual samples of the state space. The task considered in this paper is the synthesis of a polynomial controller, denoted $\pi$, that is structured as
\begin{equation}\label{eqn::controller}
\pi (x) = K(x) Z(x),
\end{equation}
where $Z \in \mathcal{P}_{d_Z}^{p}$ is a vector of monomials that are determined before optimization, and $K \in \mathcal{P}^{m \times p}_{d_K}$ is a polynomial matrix for which the coefficients are decision variables.
\begin{assumption}\label{ass:z0}
It is assumed that $Z(0)=0$.
\end{assumption}
The criteria used to determine the optimal policy, $\pi^\star$, are the `imitation' and `regularization' cost functions $\ell : \mathbb{R}^m \times \mathbb{R}^m \mapsto \mathbb{R}$ and $r: \mathcal{P}^m_{d_K + d_Z} \mapsto \mathbb{R}$, so that
\begin{equation}\label{eqn::optimal_controller}
\pi^\star := \argmin_\pi \frac{1}{N} \sum_{i=1}^N \ell (\pi(\hat{x}_i), \hat{u}_i) + r(\pi)
\end{equation}

\begin{assumption}\label{ass:convex}
The function $\ell()$ is convex in its first argument, and $r()$ is a convex function of the coefficients of $\pi$.
\end{assumption}

\subsection{Stabilizing Constraint}

Under Assumption \ref{ass:convex}, \eqref{eqn::optimal_controller} is a convex optimization problem. The extension to \eqref{eqn::optimal_controller} considered in this paper is the task of ensuring that the solution is a stabilizing w.r.t. a known polynomial system with dynamics
\begin{equation}\label{eqn::dynamics}
\dot{x} := A(x) Z(x) + B(x) u,
\end{equation}
where $A\in \mathcal{P}^{n \times n}_{d_A}$, and $B \in \mathcal{P}^{n \times m}_{d_B}$ (note that the independent variable, $t$, has been omitted for conciseness). Under closed-loop control (i.e. $u = \pi(x)$) with a controller of the form of \eqref{eqn::controller}, \eqref{eqn::dynamics} becomes\footnote{It is assumed that any discrete time artefacts (e.g. a zero-order hold in an embedded controller) are updated sufficiently quickly that \eqref{eqn::controller_dynamics} is an accurate approximation of the true system.}
\begin{equation}\label{eqn::controller_dynamics}
\dot{x} := [ A(x) + B(x) K(x) ] Z (x).
\end{equation}
Note that the LTI formulation considered in \cite{Havens2021} is a particular case of \eqref{eqn::controller_dynamics}, and occurs when $Z(x) = x$ and $d_A = d_B = d_K = 0$.

Now define $\mathcal{J} \subset \{1, \dots, n\}$ as the indices of the rows of $B(x)$ that only contain zeros, the variable $\tilde{x}$ as the $\mathcal{J}$ elements of $x$, and $P(\tilde{x} ) \in \mathcal{P}^{p \times p}_{d_P}$. For now, assume that $P(\tilde{x} )$ is positive definite $\forall \tilde{x}$, and decompose the matrix $K(x)$ in \eqref{eqn::controller} into 
$$
K(x) =: F(x) P^{-1}(\tilde{x} ),
$$
where $F \in \mathcal{P}^{m \times p}_{d_F}$. Also define $M(x)$ as the Jacobian matrix of $Z(x)$. Under the above definitions, the function $V : \mathbb{R}^n \mapsto \mathbb{R}$ defined by
$$
V(x) := Z^\top (x) P^{-1}(\tilde{x} ) Z(x)
$$
has the time derivative
\begin{multline*}
\dot{V}(x) =  Z^\top (x) \Bigg[ \sum_{j \in \mathcal{J}} \frac{\partial P^{-1}}{\partial x_j} (\tilde{x}) [A_j(x) Z(x)] \\
+ [A(x) + B(x) F(x) P^{-1}(\tilde{x})]^\top M^\top (x) P^{-1}(\tilde{x}) \\
+ P^{-1}(\tilde{x} ) M(x) [A(x) + B(x) F(x) P^{-1}(\tilde{x} ) \Bigg] Z(x),
\end{multline*}
where $x_j$ is the $j$-th element of $x$, and $A_j(x)$ is the $j$-th row of $A(x)$. Consequently, under Assumption \ref{ass:z0} the conditions
\begin{subequations}\label{eqn::PSD_conditions}
\begin{align}
P(x) & \succ 0 \label{subeqn::PSD_condition_1}\\
P(\tilde{x})A^\top (x) M^\top (x) + M(x) A(x) P(\tilde{x} ) \hspace{35pt} & \notag \\
+ F^\top (x) B^\top (x) M^\top (x) + M(x) B(x) F(x) \hspace{5pt} & \notag \\
- \sum_{j \in \mathcal{J} } \frac{\partial P}{\partial x_j} (\tilde{x}) [A_j (x) Z(x) ] & \prec 0
\end{align}
\end{subequations}
ensure that $V(x)$ is a \emph{Lyapunov function} \cite[\S 4]{khalil_2014} for the system \eqref{eqn::controller_dynamics}, which therefore implies that $u = F(x)P^{-1}(\tilde{x})Z(x)$ is a stabilizing control law for the system \eqref{eqn::dynamics}. Additionally, if $P(\tilde{x} )$ is a constant matrix, then the control law is globally stabilizing. See \cite{Prajna2004} for a more complete derivation of the above.

\subsection{Sum of Squares Approximation}

The conditions \eqref{eqn::PSD_conditions} are jointly convex in $P(x)$ and $F(x)$, as for $(P_1(x), F_1(x))$ and $(P_2(x), F_2(x))$ that both satisfy \eqref{eqn::PSD_conditions}, $(\lambda P_1(x) + (1 - \lambda ) P_2(x) , \lambda F_1(x) + (1 - \lambda ) F_2(x))$ also satisfies \eqref{eqn::PSD_conditions} for all $\lambda \in [0,1]$. Despite this desirable property, \eqref{eqn::PSD_conditions} remains computationally intractable: for example, in the case where $p=1$, condition \eqref{subeqn::PSD_condition_1} equates to a polynomial nonnegativity constraint. \emph{Sum of squares} techniques are commonly used to render polynomial nonnegativity constraints tractable\footnote{See \cite{Parrilo} for a more complete description of sum of squares approaches and the challenges of polynomial nonnegativity constraints; only sufficient detail to understand the following is presented here.}: a polynomial $f \in \mathcal{P}$ is said to be sum of squares iff $f(x) = \sum_{i=1}^{\hat{N}} [f_i(x)]^2$, where each $f_i(x) \in \mathcal{P}$ for all $i \in \{0, \dots, \hat{N} \}$. Equivalently, a polynomial of degree $2d$ is sum of squares if and only if there exists a vector of monomials up to degree $d$, denoted $z(x)$, and a positive semidefinite matrix, $Q$, such that $p(x) = z^\top (x) Q z(x)$. Clearly, if a polynomial is sum of squares then it is also positive semidefinite. The following useful result is proven in \cite{Prajna2004}:

\begin{lemma}
Given a $p \times p$ polynomial matrix, $\hat{P}(x)$, and a $p$-dimensional vector of monomials, $\hat{Z}(x)$, if there exists $Q \succeq 0$ such that
$$
v^\top \hat{P}(x) v = [\hat{Z}(x) \otimes v ]^\top Q [\hat{Z}(x) \otimes v]
$$
for $v \in \mathbb{R}^p$, then $\hat{P}(x) \succeq 0 $ for all $x$.
\end{lemma}
Consequently, the system \eqref{eqn::controller_dynamics} is stable w.r.t. origin if there exists vectors of monomials\footnote{Note that the constituent monomials of $z_1(x)$ and $z_2(x)$ must be chosen to match the degree of the matrix polynomials on the L.H.S. of \eqref{eqn::constraint_1} and \eqref{eqn::constraint_2}}, $z_1(x)$ and $z_2(x)$, and matrices $Q_1$ and $Q_2$ such that
\begin{align}
v^\top [P(\tilde{x} ) - \epsilon_1 I] v &= [z_1(x) \otimes v]^\top Q_1 [z_1(x) \otimes v] \label{eqn::constraint_1} \\
v^\top \Bigg[ P(\tilde{x})A^\top (x) M^\top (x) & + M(x) A(x) P(\tilde{x} ) \notag  \\
 + F^\top (x) B^\top (x)  M^\top &  (x) + M(x) B(x) F(x) \notag \\
  - \sum_{j \in \mathcal{J} }  \frac{\partial P}{\partial x_j} & (\tilde{x})  [A_j (x) Z(x) ] + \epsilon_2 I \Bigg] v \notag \\
   & = [z_2(x) \otimes v]^\top Q_2 [z_2(x) \otimes v]  \label{eqn::constraint_2} \\
Q_1 & \succeq 0  \label{eqn::constraint_3} \\
Q_2 & \preceq 0  \label{eqn::constraint_4}
\end{align}
where $v \in \mathbb{R}^p$ and $\epsilon_1, \epsilon_2 > 0$. The expressions on either side of the equalities in \eqref{eqn::constraint_1} and \eqref{eqn::constraint_2} are polynomials in terms of $x$ and $v$, and the equalities therefore have the effect of constraining the coefficients of each constituent monomial to be equal on both sides. Constraints of this form are not readily parsed by general purpose optimization software, so an equivalent representation in terms of coefficient equalities is presented in Appendix \ref{section::monomial_representation} (due to space constraints, only \eqref{eqn::constraint_1} is addressed explicitly, but \eqref{eqn::constraint_2} follows a very similar approach).

\subsection{Stabilizing Imitation Learning}

Imitation learning of a control law \eqref{eqn::controller} that is stabilizing for \eqref{eqn::controller_dynamics} can therefore be expressed as the optimization problem
\begin{equation}\label{eqn::stabilizing_imitation_learning}
\begin{aligned}
\min_{\{ K_i \}, \{F_i \}, \{P_i \}, Q_1, Q_2 }\  & \frac{1}{N} \sum_{i=1}^N \ell (K(\hat{x}_i ) Z(\hat{x}_i), \hat{u}_i) + r( \{ K_i \} ) \\
\textrm{s.t.} \ & K(x) P(\tilde{x} ) = F(x) \\
& \eqref{eqn::constraint_1} - \eqref{eqn::constraint_4},
\end{aligned}
\end{equation}
where $\{ P_i \}$, $\{ F_i \}$, and $\{K_i \}$ are used as shorthand notation for the matrices $P_i \in \mathbb{S}^{p}$, $F_i \in \mathbb{R}^{m \times p}$, and $K_i \in \mathbb{R}^{m \times p}$ in the monomial representations
\begin{gather*}
P(\tilde{x}) = \sum_{i=1}^{|\mathcal{M}_{d_P}[\tilde{x}] |} P_i  \mathcal{M}_{d_P}[\tilde{x}]_i, \\
F(x) = \hspace{-5pt} \sum_{i=1}^{|\mathcal{M}_{d_F}[x]|}  \hspace{-5pt} F_i \mathcal{M}_{d_F}[x]_i, \quad
K(x) = \hspace{-5pt} \sum_{i=1}^{|\mathcal{M}_{d_K}[x]|}  \hspace{-5pt} K_i \mathcal{M}_{d_K}[x]_i.
\end{gather*}
Under Assumption \ref{ass:convex}, the objective function of \eqref{eqn::stabilizing_imitation_learning} is convex in $\{K_i \}$, and the constraints \eqref{eqn::constraint_1} - \eqref{eqn::constraint_4} are convex in all decision variables, but the constraint $K(x)P(\tilde{x} ) = F(x)$ is \emph{biaffine}, so problem \eqref{eqn::stabilizing_imitation_learning} is \emph{biconvex} in $\{K_i\}$ and $(\{F_i\}, \{P_i\}, Q_1, Q_2)$.

\section{Optimization Algorithms}

The `Kalman constraint' in \cite{Palan2020} is also biaffine, and ADMM was proposed as a solution heuristic in that paper as it has been shown to be convergent (to a not necessarily optimal point) for this problem class when using a sufficiently large penalty parameter \cite{Gao2020}. A projected gradient descent algorithm was proposed as an alternative to ADMM in \cite{Havens2021}, where it was shown to provide superior performance in some conditions. This section therefore presents both an ADMM and a projected gradient descent algorithm as heuristics for the solution of \eqref{eqn::stabilizing_imitation_learning}.

\subsection{Alternating Direction Method of Multipliers}
The constraint $K(x) P(\tilde{x} ) = F(x) $ in \eqref{eqn::stabilizing_imitation_learning} is equivalent to
\begin{multline*}
 \hspace{-5pt} \sum_{i=1}^{|\mathcal{M}_{d_K}[x]|}  \sum_{j=1}^{|\mathcal{M}_{d_P}[\tilde{x}] |}     K_i P_j \mathcal{M}_{d_K}[x]_i \mathcal{M}_{d_P}[\tilde{x}]_j  \\
 =  \sum_{k=1}^{|\mathcal{M}_{d_F}[x]|}   F_i \mathcal{M}_{d_F}[x]_k,
\end{multline*}
which is in turn equivalent to $| \mathcal{M}_{d_F} [x] |$ constraints of the form 
$$
\sum_{(i, j) \in \mathcal{E}[k] } K_i P_j = F_k
$$
where $\mathcal{E}[k] := \{ (i, j) : \mathcal{M}_{d_K}[x]_i \mathcal{M}_{d_P}[\tilde{x}]_j = \mathcal{M}_{d_F}[x]_k  \}$. Therefore, the scaled augmented Lagrangian for \eqref{eqn::stabilizing_imitation_learning} is defined as 
\begin{multline*}
\mathcal{L}_\rho (\{K_i \}, \{F_i \}, \{P_i \}, Q_1, Q_2 ) \\ 
:= J(\{ K_i \} ) + \sum_{k=1}^{ | \mathcal{M}_{d_F} [x] | } \frac{\rho}{2} \left\| F_k - \sum_{(i, j) \in \mathcal{E}[k] } K_i P_j + Y_k \right\|_2^2,
\end{multline*}
where $Y_k \in \mathbb{R}^{m \times p}$ are scaled dual variables defined for all $k \in \{0, \dots, |\mathcal{M}_{d_F} [x] | \}$, and
$J(\{K_i \} ) := \frac{1}{N} \sum_{i=1}^N \ell (K(\hat{x}_i ) Z(\hat{x}_i), \hat{u}_i) + r( \{ K_i \} )$. Given initial values $(\{ K^{(0)}_i \}\{F^{(0)}_i \}, \{P^{(0)}_i \}, Q_1^{(0)}, Q_2^{(0)} )$, the ADMM algorithm is then
\begin{align*}
\{ K^{(l+1)} \} & := \argmin_{\{ K_i \}} \mathcal{L}_\rho (\{K_i \}, \{F_i^{(l)} \}, \{P_i^{(l)} \}, Q_1^{(l)}, Q_2^{(l)} ),  \\
& \hspace{-42pt} \left( \{F_i^{(l+1)} \}, \{P_i^{(l+1)} \}, Q_1^{(l+1)}, Q_2^{(l+1)} \right) \\
& := \hspace{-10pt} \argmin_{\{F_i \}, \{P_i \}, Q_1, Q_2} \hspace{-5pt} \mathcal{L}_\rho (\{K_i^{(l+1)} \}, \{F_i \}, \{P_i \}, Q_1, Q_2 ) \\
& \hspace{40pt} \mathrm{s.t.} \ \eqref{eqn::constraint_1} - \eqref{eqn::constraint_4}, \\
Y_k^{(l+1)} & := Y_k^{(l)} +  F_k^{(l+1)} - \sum_{(i, j) \in \mathcal{E}[k] } K_i^{(l+1)} P_j^{(l+1)} \quad \forall k.
\end{align*}
\subsection{Projected Gradient Descent}\label{section::grad_desc}
Problem \eqref{eqn::stabilizing_imitation_learning} is equivalent to 
\begin{equation}\label{eqn::stabilizing_imitation_learning_2}
\begin{aligned}
\min_{\{ K_i \}, \{F_i \}, \{P_i \}, Q_1, Q_2 }\  & \frac{1}{N} \sum_{i=1}^N \ell (F(\hat{x}_i) P^{-1}(\tilde{\hat{x}}_i ) Z(\hat{x}_i), \hat{u}_i) \\ & \hspace{60pt} + r( m(\{ F_i \}, \{ P_i \} ) ) \\
\textrm{s.t.} \ & \eqref{eqn::constraint_1} - \eqref{eqn::constraint_4},
\end{aligned}
\end{equation}
where $m$ is a nonlinear function that maps the values of $\{ F_i \}$ and $\{ P_i \}$ to $\{K_i \}$ such that $K(x)  = F(x) P^{-1}(\tilde{x} )$. It is assumed that the objective function in \eqref{eqn::stabilizing_imitation_learning_2} is differentiable almost everywhere, so that, given initial values $(\{ K^{(0)}_i \}\{F^{(0)}_i \}, \{P^{(0)}_i \}, Q_1^{(0)}, Q_2^{(0)} )$, a solution to \eqref{eqn::stabilizing_imitation_learning_2} can be approximated using the projected gradient descent algorithm:
\begin{align*}
\quad \tilde{F}^{(l)}_j & := F_j^{(l)} - \alpha \frac{\partial J}{\partial F_j} \left( \{F^{(l)}_i \}, \{P^{(l)}_i \} \right) \quad \forall j, \\
\tilde{P}^{(l)}_j & := P_j^{(l)} - \alpha \frac{\partial J}{\partial P_j} \left( \{F^{(l)}_i \}, \{P^{(l)}_i \} \right) \quad \forall j,  \\
& \hspace{-30pt} \left(\{F_i^{(l+1)} \}, \{P_i^{(l+1)} \}, Q_1^{(l+1)}, Q_2^{(l+1)} \right)\\
& := \argmin_{\{F_i\}, \{P_i \}} \ \sum_{j=1}^{|\mathcal{M}_{d_F}[x]|} \left\| \tilde{F}^{(l)}_j - F_j \right\|_F^2 \\ 
& \hspace{100pt} + \sum_{j=1}^{|\mathcal{M}_{d_P}[x]|} \left\| \tilde{P}_j^{(l)} - P_j \right\|_F^2\\
& \hspace{38pt} \mathrm{s.t} \  \eqref{eqn::constraint_1} - \eqref{eqn::constraint_4},
\end{align*}
where $\alpha$ is a step-size parameter, and $\frac{\partial J}{\partial \gamma} \left( \{F^{(i)} \}, \{P^{(i)} \} \right)$ is the gradient of the objective function of \eqref{eqn::stabilizing_imitation_learning_2} w.r.t. given variable $\gamma$ evaluated at $ \left( \{F^{(i)} \}, \{P^{(i)} \} \right)$. The gradients could also be evaluated w.r.t. a random subsample of the data \eqref{eqn::data}, producing a stochastic approximation (for which the expected value is equal to the deterministic gradient) that can be useful for helping gradient-based optimization methods escape from local minima.

\section{Numerical Experiments}

This section demonstrates the safe imitation learning approach proposed in this paper on two systems: a nonlinear system with a linear feedback controller, and a linear system with a nonlinear feedback controller. For both sets of experiments, training data was generated for $i \in \{1, \dots, N \}$ by randomly sampling $\hat{x}_i \in [-10, 10]^n$, then computing the expert controller with a random perturbation as $\hat{u}_i = K(\hat{x}_i) Z(\hat{x}_i ) + \epsilon$, where $\epsilon$ was normally distributed with mean 0 and covariance matrix $\sigma I$ ($\sigma =1$ was used for all experiments). The least square cost was used for the imitation loss function $\ell()$, and the regularization loss was not used. Both the ADMM algorithm and projected gradient descent algorithm were run for a fixed number of iterations, and the parameters $\rho$ and $\alpha$ were tuned manually. For each system, every element of $\{K_i\}$, $\{P_i\}$, and $\{F_i\}$ was initialized by sampling uniformly from $[-5, 5]$, and the optimization process was completed using ten random number seeds for each of $N \in \{10, 100, 1000\}$.

The experiments were implemented in Python 3.8 on an 2.60GHz Intel i7-9750 CPU. The polynomial constraints were parsed (as described in Appendix \ref{section::monomial_representation}) using Sympy \cite{Sympy}. All optimization subproblems for both the projected gradient descent algorithm and ADMM were modelled using CVXPY \cite{diamond2016cvxpy} and solved with SCS \cite{scs}. The derivatives for the projected gradient descent algorithm were calculated using Jax \cite{jax2018github}. The code for the experiments is publicly available at \texttt{github.com/sebastian-east/sos-imitation\\-learning}.

\subsection{Nonlinear System}
\begin{figure}
	\centering
	\input{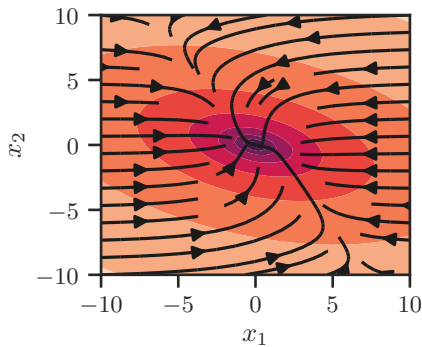}
\caption{Contour plot of the learned Lyapunov function and streamlines of the closed-loop system using the learned controller from a single ADMM experiment. \label{fig::lyapunov1}}
\end{figure}

The first system considered was taken from \cite{Prajna2004}, for which $n=2$, $m=1$, and the dynamics are given by
\begin{gather*}
A(x) = \begin{bmatrix}
-1 + x_1 - \frac{3}{2} x_1^2 - \frac{3}{4} x_2^2 & \frac{1}{4} - x_1^2 - \frac{1}{2} x_2^2 \\ 0 & 0
\end{bmatrix} \\ 
Z(x) = [x_1, x_2 ]^\top, \quad B(x) = [0, 1]^\top.
\end{gather*}
The `expert' took the form of a linear state-feedback controller with $K(x) = [-2, -10]$, which can be shown to be stabilizing using a second order sum of squares Lyapunov function. Figure \ref{fig::training1} shows the imitation loss during the training process for both proposed optimization algorithms, using parameters $\rho=1$ and $\alpha = 10^{-5}$. The parameters of the learning controller were set at $d_F = 0$ and $d_P = 0$ (i.e. the same as the expert). It is clearly shown that ADMM converges to a close approximation of the expert controller within a few iterations, across all seeds and values of $N$, whereas the projected gradient descent algorithm often gets stuck in local minima, and requires significantly more iterations even when it does converge to a low imitation loss. The performance of projected gradient descent shown in Figure \ref{fig::training1} is generally representative of the best performance found for a broad range of tested values of $\alpha$. Figure \ref{fig::lyapunov1} shows the contours of the Lypunov function learned using ADMM for a single example, and clearly shows that its value decreases along the accompanying closed-loop system trajectories.

\begin{figure}
	\centering
	\input{training1.pgf}
\caption{Imitation loss at each iteration of the learning process for the nonlinear system using ADMM and projected gradient descent. Each line represents the results for a single random initialization. \label{fig::training1}}
\end{figure}

\subsection{Nonlinear Control}

The second system considered had the (marginally stable) linear dynamics
\begin{gather*}
A(x) = \begin{bmatrix}
0 & 1\\ -1  & 0
\end{bmatrix}, \quad Z(x) = [x_1, x_2 ]^\top, \quad B(x) = [0, 1]^\top,
\end{gather*}
for which the `expert' took the form of a 3rd order state-feedback controller with $K(x) = [-0.1 - 0.1 x_1^2, -0.1 - 0.1x_2^2]$, which can be shown to be stabilizing using a fourth order sum of squares Lyapunov function. The learning parameters were set this time at $\rho = 1000$ and $\alpha = 10^{-8}$ (the values used for the previous experiment were found to perform poorly), and the parameters of the learning controller were set at $d_F = 2$ and $d_P = 0$. Note that these controller parameters only allow a learned Lyapunov function of degree two, and that no second order Lyapunov function exists for this system when controlled by the expert, so it is impossible for the algorithms to converge on the `correct' controller. This limitation could be resolved using an alternative choice of $Z(x)$, but was kept to illustrate the learning process in the presence of a `correspondence' issue. 

Figure \ref{fig::training2} shows the imitation loss during training for both algorithms. The ADMM algorithm now required up to 200 iterations to converge, and this time converged to a significantly higher imitation loss than for the previous experiments. This is possibly due to the aforementioned correspondence issue, but also may be caused by the fact that the nonconvex constraint $K(x) P(\tilde{x}) = F(x)$ now contains more terms, which may be more challenging to address with ADMM. The comparative performance of the projected gradient descent algorithm was again worse, requiring a much higher number of iterations, and often immediately becoming entrenched in a local minimum with high imitation loss. The illustrated value of $\alpha$ is again representative of the best performance that could be found across a wide range of values, but performance could possibly be improved by updating $\alpha$ during the learning process, or by using stochastic gradient descent (as discussed in Section \ref{section::grad_desc}).

\begin{figure}
	\centering
	\input{training2.pgf}
\caption{Imitation loss at each iteration of the learning process for the nonlinear system using ADMM and projected gradient descent. Each line represents the results for a single random initialization. \label{fig::training2}}
\end{figure}
\section{Conclusion}

This paper demonstrated an imitation learning method that can guarantee the stability of a learned polynomial policy when used to control a polynomial system. An ADMM algorithm and projected gradient descent algorithm were proposed as heuristic solutions to the associated optimization problem, and it was demonstrated through numerical experiments that ADMM is generally more effective.

One of the limitations of the proposed approach is that it generates a globally stabilizing controller: global stabilization may be conservative, and it was found that the constraints \eqref{eqn::constraint_1}-\eqref{eqn::constraint_4} are often infeasible for more complex systems than those presented here. Future work will therefore investigate the extension of this approach to consider stability on bounded domains.

\appendix

\subsection{Monomial Representation of Polynomial Equality Constraints}\label{section::monomial_representation}

To aid the following, consider the vector-valued polynomials $v_1(x) \in \mathcal{P}_{d_1}^p$ and $v_2(x) \in \mathcal{P}_{d_2}^p$, and the vector valued matrix $H(x) \in \mathcal{P}_{d_3}^{p \times p}$. The matrices $H(x)$ and $v_2(x) v_1^\top(x)$ can be decomposed into the monomial representations
\begin{gather*}
H(x) = \sum_{j=1}^{|\mathcal{M}_{d_3}[x] |} H_j  \mathcal{M}_{d_3}[x]_j \\
v_2(x) v_1^\top (x) = \sum_{i=1}^{|\mathcal{M}_{d_1 + d_2}[x] |} C_i \mathcal{M}_{d1 + d2}[x]_i,
\end{gather*}
where $H_j,C_i \in \mathbb{R}^{p \times p}$ $\forall \ i,j$. Therefore, it can be shown that the polynomial $v_1^\top(x)  H(x) v_2(x)$ has the equivalent monomial representation
\begin{multline}\label{eqn::monomial_representation}
v_1^\top (x) H (x) v_2(x) = \tr (v_2(x) v_1^\top (x) H(x)) \\ 
=  \sum_{i=1}^{|\mathcal{M}_{d_1 + d_2}[x] |} \sum_{j=1}^{|\mathcal{M}_{d_3}[x] |} \tr ( C_i^\top H_j) \mathcal{M}_{d1 + d2}[x]_i \mathcal{M}_{d_3}[x]_j.
\end{multline}

Now consider \eqref{eqn::constraint_1}. The following terms have the monomial representations
\begin{align*}
\quad vv^\top &= \sum_{i=1}^{|\mathcal{M}_2[v]|} C_i^{(1)} \mathcal{M}_2[v]_i, \\
[z_1(x) \otimes v] [z_1(x) \otimes v]^\top &= \hspace{-5pt} \sum_{k=1}^{|\mathcal{M}_2[z_1(x) \otimes v] |} \hspace{-5pt} C^{(2)}_k  \mathcal{M}_2[z_1(x) \otimes v]_k,
\end{align*}
where $C_i^{(1)} \in \mathbb{S}^{p} \ \forall i$ and $C_k^{(2)} \in \mathbb{S}^{p^2 }\ \forall k$ are constants. Therefore, using the identity \eqref{eqn::monomial_representation}, the constraint \eqref{eqn::constraint_1} is equivalent to 
\begin{multline}
\sum_{i=1}^{|\mathcal{M}_2[v]|} \sum_{j=1}^{|\mathcal{M}_{d_P}[x]|} \tr (C_i^{(1)} P_j ) \mathcal{M}_2[v]_i \mathcal{M}_{d_P}[x]_j - \sum_{i=1}^n \epsilon_1 v_i^2 \\
= \sum_{i=1}^{|\mathcal{M}_2[z_1(x) \otimes v]|} \tr (C_k^{(2)} Q) \mathcal{M}_2[z_1(x) \otimes v]_k.
\end{multline}
Hence, \eqref{eqn::constraint_1} can be represented by $|\mathcal{M}_2[z_1(x) \otimes v]|$ equality constraints of the form
$$
\tr (C_i^{(1)} P_j ) - \epsilon_1 \mathcal{I}_{v^2}(i) = \tr (C_k^{(2)} Q ),
$$
where
$$
\mathcal{I}_{v^2} (i) : =  \begin{cases} 1 & \mathcal{M}_2[v]_i \in \{v_1^2, \dots, v_n^2 \} \\
0 & \textrm{otherwise} \end{cases},
$$
and $k$ is implicitly defined for each pair $i,j$ such that $
\mathcal{M}_2[v]_i \mathcal{M}_{d_P}[x]_j = \mathcal{M}_2[z_1(x) \otimes v]_k$.

A similar approach can be used for the constraint \eqref{eqn::constraint_2} by noting that 
$$
\frac{\partial P }{\partial x_j} (\tilde{x}) = \sum_{j=1}^{|\mathcal{M}_{d_P}[\tilde{x}] |} P_j  \frac{\partial \mathcal{M}_{d_P}[\tilde{x}]_j} {\partial x_j}
$$
for all $j \in \mathcal{J}$.

\bibliographystyle{ieeetr}
\bibliography{bibliography}

\end{document}